\documentclass[12pt]{amsart}

\usepackage{verbatim}
\usepackage{amssymb}
\usepackage{upref}
\usepackage[all]{xy}
\usepackage{tikz,pgf}

\newtheorem{thm}{Theorem}[section]

\newtheorem{prop}[thm]{Proposition}

\theoremstyle{definition}

\newtheorem{ex}[thm]{Example}

\newtheorem{rem}[thm]{Remark}

\numberwithin{equation}{section}

\newcommand{\thmref}[1]{Theorem~\textup{\ref{#1}}}

\newcommand{\midtext}[1]{\quad\text{#1}\quad}

\renewcommand{\and}{\midtext{and}}

\newcommand{\C}{\mathbb C}

\newcommand{\CC}{\mathcal C}
\newcommand{\KK}{\mathcal K}

\newcommand{\GG}{\mathcal G}
\newcommand{\LL}{\mathcal L}

\newcommand{\OO}{\mathcal O}

\newcommand{\p}{\phi}

\newcommand{\x}{\xi}

\newcommand{\Chi}{\raisebox{2pt}{\ensuremath{\chi}}}
\renewcommand{\epsilon}{\varepsilon}

\newcommand{\G}{\Gamma}

\DeclareMathOperator{\obj}{Obj}

\DeclareMathOperator*{\spn}{span}
\DeclareMathOperator*{\clspn}{\overline{\spn}}

\newcommand{\case}{& \text{if }}

\newcommand{\<}{\langle}
\renewcommand{\>}{\rangle}

\renewcommand{\bar}{\overline}

\let\oldmarginpar\marginpar        
\newcommand\sknote[1]%
{\textbf{*}\-\oldmarginpar[\raggedleft\tiny #1 --sk]%
{\raggedright\tiny #1 --sk}}

\newcommand{\mdot}{\cdot}

\begin{document}

\title{Characterizing graph $C^*$-correspondences}

\author{S.~Kaliszewski}
\author{Nura Patani}
\author{John Quigg}



\date{\today}

\begin{abstract}
Every separable nondegenerate $C^*$-cor\-re\-spon\-dence 
over a commutative $C^*$-algebra with discrete spectrum 
is isomorphic to a graph correspondence.
\end{abstract}

\maketitle

Let $E$ be a directed graph with vertex set $V$,
edge set $E^1$, and range and source maps $r,s:E^1\to V$.
The \emph{graph correspondence} is the nondegenerate
$C^*$-correspondence  
$X_E$ over $c_0(V)$ defined by
\[
X_E=\bigl\{\x\colon E^1\to\C \bigm|
\text{the map }v\mapsto \sum_{s(e)=v}|\x(e)|^2
\text{ is in }c_0(V)
\bigr\},
\]
with module actions and $c_0(V)$-valued inner product given by
\[
(a\mdot\x\mdot b)(e)=a(r(e))\,\x(e)\,b(s(e))
\quad\text{and}\quad
\<\x,\eta\>(v)=\sum_{s(e)=v}\bar{\x(e)}\eta(e)
\]
for $a,b\in c_0(V)$,  $\x,\eta\in X_E$, $e\in E^1$, and $v\in V$.
(We cite \cite{raeburngraph} as a general reference on graph algebras, and \cite[Chapter~8]{raeburngraph} in particular for graph cor\-re\-spon\-dences.)

A graph correspondence $X_E$ contains 
at least as much information about $E$ as the graph algebra $C^*(E)$, 
since the Cuntz-Pimsner algebra $\OO_{X_E}$ 
is isomorphic to $C^*(E)$ 
(\cite[Example~8.13]{raeburngraph};
see also \cite[Example 1, p.~4303]{katsgen}).  
Indeed, many properties of $E$ 
are directly reflected in properties of $X_E$:
for example, 
$X_E$ is \emph{full} in the sense that $\clspn\<X_E,X_E\>=c_0(V)$ 
if and only if the graph $E$ has no sinks;
the homomorphism $c_0(V)\to \LL_{}(X_E)$ 
associated to the left module operation 
maps into the compacts $\KK_{}(X_E)$
if and only if no vertex receives infinitely many edges;
$c_0(V)$ maps faithfully into $\LL_{}(X_E)$
if and only if $E$ has no sources.

In this paper we further investigate 
this connection between graphs and $C^*$-correspondences.  
In Section~\ref{characterization}, 
our main result (Theorem~\ref{characterization})
expands and elaborates on a remark of Schweizer 
(\cite[\S1.6]{schcrosspim}; see also \cite[Chapter~1, Example~2]{pim})
to the effect that every $C^*$-correspondence over $C^*$-algebras $c_0(X)$ and $c_0(Y)$ is unitarily equivalent to a correspondence arising from a ``diagram'' from $X$ to $Y$. 
(Schweizer's diagrams can be put in our context by taking $V=X\cup Y$.)
Specifically, 
given a (separable, nondegenerate) $C^*$-correspondence $X$ over $c_0(V)$, 
we construct a graph $E$ 
such that $X_E\cong X$ as $C^*$-correspondences.
We also show (Proposition~\ref{recover-F}) that 
for any graph $F$ with vertex set $V$, 
applying our construction to $X_F$ gives a graph $E$ such that 
$E\cong F$ via a vertex-fixing graph isomorphism.

In Section~\ref{functoriality} we show that the assignment $E\mapsto X_E$
can be extended to a functor $\G$ between certain categories of graphs and correspondences.  $\G$ is very nearly a category equivalence: it is
essentially surjective, faithful, and ``essentially full''
(see Theorem~\ref{close-to-equiv}), but not full.  
It is also injective on objects, and reflects isomorphisms
(Theorem~\ref{injective}).

Although our construction of a graph from a $C^*$-correspondence is fairly elementary, we believe our results give valuable insight into graph correspondences.  
In work currently in progress, we are applying 
the techniques illustrated here to the more general context of topological graphs.
We are also interested in the analogous problem for $k$-graphs,
although since there are $k$-graph $C^*$-algebras 
which are not isomorphic to any directed-graph $C^*$-algebra,
our main result in this paper implies that
for $k$-graphs it will not be adequate to consider $C^*$-correspondences 
over $c_0$ of the vertex set.

\subsection*{Notation and Terminology}

Throughout this paper, $V$ will denote a countable set.
All graphs considered will have vertex set $V$,
and (by definition) countably many edges.
We use $\Chi_S$ to denote the characteristic function of a set $S$,
and we abbreviate this to $\Chi_x$ when $S=\{x\}$ is a singleton.
We will write $p_v$ for the minimal projection $\Chi_v\in c_0(V)$.

A $C^*$-correspondence $X$ over a $C^*$-algebra $A$ is said to be \emph{nondegenerate}
if the set of products $S= \{ a\mdot \xi \mid a\in A, \xi\in X \}$ has dense span in $X$. 

If $X$ and $Y$ are $C^*$-correspondences over $A$, 
a morphism from $X$ to $Y$ is a  linear map
$\psi\colon X\to Y$ satisfying
\begin{equation}\label{psi}
\psi(a\mdot\x\mdot b)=a\cdot\psi(\x)\cdot b\and
\<\psi(\x),\psi(\eta)\>=\<\x,\eta\>
\end{equation}
for $a,b\in A$ and $\x,\eta\in X$.  
Note that such morphisms will be isometric, and in particular injective.

Note that for a graph $E$,
$X_E$ is densely spanned by the linearly independent set 
$\{\chi_e \mid {e \in E^1} \}$, 
whose elements we refer to as \emph{generators} of the correspondence $X_E$.
The $c_0(V)$-valued inner product on $X_E$ is characterized on generators by
\[
\<\chi_e,\chi_f\>=
\begin{cases}
p_{s(e)}\case e=f\\
0\case e\ne f.
\end{cases}
\]

\section{The Isomorphisms}
\label{characterization}

\begin{thm}
\label{characterize}
Let $V$ be a countable set.
Every separable nondegenerate $C^*$-correspondence over $c_0(V)$ 
is isomorphic to 
the graph correspondence of a directed graph with vertex set $V$.
\end{thm}

\begin{proof}
For ease of notation, set $A=c_0(V)$.
Let $X$ be a separable nondegenerate $C^*$-correspondence over $A$, 
For each $v\in V$ put
\[
X_v=X\mdot p_v.
\]
Because finite sums of the $p_v$'s form a bounded approximate
identity in $A$ 
(and the right action of $A$ on $X$ is automatically nondegenerate),
each $\xi\in X$ can be approximated by sums of the form
\[
 \xi \mdot\Bigl(\sum_{v\in F} p_v \Bigr)
= \sum_{v\in F}\xi\mdot p_v  \in \spn \bigcup_{v\in F} X_v,
\]
with $F\subseteq V$ finite; 
so the $X_v$'s span a dense subspace of $X$.

Moreover,
for  $u\ne v$ in $V$ we have $p_u p_v=0$ in $A$,
and hence
for $\x\in X_u$ and $\eta\in X_v$
we have
\[
\<\x,\eta\>
=\<\x\mdot p_u,\eta\mdot p_v\>
=p_u\<\x,\eta\>p_v
=\<\x,\eta\>p_up_v
= 0.
\]
It follows by a standard argument that the subspaces $X_v$ are linearly independent.

Now with the norm inherited from $X$, 
each $X_v$ is in fact a Hilbert space,
with inner product $\langle\cdot,\cdot\rangle_v$
(conjugate-linear in the first variable)
given by
\[
\langle\xi,\eta\rangle_{v}=\<\xi,\eta\>(v).
\]
For $u\in V$ put
\[
 X_{uv} = p_u\mdot X\mdot p_v = p_u\mdot X_v \subseteq X_v.
\]
Then for
$u\ne w$ in $V$ 
the closed subspaces
$X_{uv}$ and $X_{wv}$ of $X_v$ are orthogonal
because, for $\x\in X_{uv}$ and $\eta\in X_{wv}$, 
\[
\<\xi,\eta\>_v
=\<\xi,\eta\>(v)
=\<p_u\cdot\xi, p_w\cdot\eta\>(v)
=\<\xi, p_u p_w\cdot\eta\>(v) = 0.
\]
Because the left action of $A$ on $X$ (hence on $X_v$) 
is nondegenerate, 
(arguing as above)
the family $\{X_{uv}\mid u\in V\}$ of mutually orthogonal subspaces
spans a dense subspace of $X_v$. 

For each $u,v\in V$ choose
an orthonormal basis
$E_{uv}$
of the
Hilbert space $X_{uv}$, and put
\[
E^1=\bigcup_{u,v\in V}E_{uv}.
\]
Since $X$ is separable,
every $E_{uv}$ is countable 
and hence $E^1$ is countable.  
Since the subspaces $X_{uv}$ are linearly independent, 
the union above is a disjoint union, so we can define $r,s:E^1\to V$ by
\[
r(e)=u\midtext{and}s(e)=v\midtext{if}e\in E_{uv}.
\]
This gives a graph $E=(V,E^1,r,s)$.
\medskip

We now show that $X\cong X_E$ as $C^*$-correspondences
over $A$.  
For this we need a linear surjection $\psi:X\to X_E$ 
satisfying~\eqref{psi}.
Note that by construction, 
$E^1$ is a linearly independent subset of $X$,
and for each $e\in E^1$ we have $\Chi_e\in X_E$ by definition of $X_E$.  
Thus, setting $X_0=\spn E^1$, 
we can define $\psi\colon X_0 \to X_E$ by linearly extending the rule
\[
 \psi(e) = \Chi_e
\]
from $E^1$ to $X_0$.

Now it is straightforward to check the action and inner product identities
\eqref{psi}
for $a=p_u$, $b=p_v$, and $\xi,\eta\in E^1$, 
and the identities then extend to $a,b\in c_c(V)$ and $\xi,\eta\in X_0$ by linearity.  In particular, $\psi$ is isometric on $X_0$, so since $X_0$ is dense in $X$, $\psi$ extends to an isometric $c_0(V)$-bimodule map on $X$ by a standard argument.  Moreover, 
$\psi$ is surjective because the image of $E^1$ has dense span in $X_E$.
\end{proof}

\begin{ex}\label{sigma}
Let $V$ be a countable set, 
and let $\sigma$ be a map of $V$ into itself.
This gives rise to a natural $C^*$-correspondence $X$ over $A=c_0(V)$
defined by $X = c_0(V)$, with actions and inner product given in terms of the usual operations in $c_0(V)$ by
\[
 a\mdot \xi\mdot b = (a\circ\sigma)\xi b
\quad\text{and}\quad
\langle \xi,\eta \rangle = \bar\xi\eta
\]
for $a,b\in A$ and $\xi,\eta\in X$.
Note that $a\circ\sigma\in \ell^\infty(V) = M(c_0(V))$, so the left action is well-defined.
(This construction motivated the more general discussion in \cite[\S1.2]{schcrosspim}.)

Following the construction and 
notation of the proof of \thmref{characterize},
we have
\[
X_v = X\mdot p_v = \{ \xi\Chi_v \mid \xi\in c_0(V) \}
= \spn \Chi_v,
\]
and then, since the left action of $p_u\in A$ on $\Chi_v\in X$ is $(\Chi_u\circ\sigma)\Chi_v 
= \Chi_{\sigma^{-1}(u)}\Chi_v
= \Chi_v$ if $\sigma(v)=u$, and is zero otherwise, we have
\[
X_{uv} = 
\begin{cases} 
\spn \Chi_v & \text{if }\sigma(v)=u \\ 
0 & \textrm{ otherwise.} \end{cases}
\]
So in this example, there is a natural choice of
orthonormal basis for each $X_{uv}$ 
(provided $\sigma(v)=u$), 
namely $E_{uv} = \{ \Chi_v \}$.
(Each $\Chi_v$ is a unit vector in $X$ because
the $C^*$-correspondence norm on $X$ 
is just the supremum norm on $c_0(V)$.)
Since $E_{uv} = \emptyset$ when $\sigma(v)\neq u$,
the resulting graph $E=(V,E^1,r,s)$
associated to $X$ has 
\[
E^1 = \{ \Chi_{v} \mid v\in V \},\quad
r(\Chi_{v}) = \sigma(v),\quad
s(\Chi_{v}) = v.
\]

Note that $E$ is isomorphic to the graph $E_\sigma$ with vertex set $V$ which is defined using $\sigma$ in a natural way by defining an edge from $v$ to $u$ precisely when $\sigma(v)=u$:
\[
E^1_\sigma = \{(\sigma(v),v)\mid v\in V \},\quad
r(\sigma(v),v)=\sigma(v),\quad
s(\sigma(v),v)=v.
\]
\end{ex}

\begin{ex}
Let $V=\{u,v\}$ and define $\sigma:V\to V$ 
by $\sigma(u) = v$ and $\sigma(v) = u$.
Let $X$ be the $C^*$-correspondence over $c_0(V)$ associated to $\sigma$
as in Example~\ref{sigma}.  
Then $X_{uu}$ and $X_{vv}$
are zero-dimensional, 
and we can choose unit vectors 
$\Chi_v \in X_{uv} = \spn \Chi_v$
and  $\Chi_u\in X_{vu}=\spn \Chi_u$.  
We have $E^1=\{\Chi_u, \Chi_v\}$
with $r(\Chi_{v}) = s(\Chi_{u}) = u$
and $s(\Chi_{v}) = r(\Chi_{v}) = v$, 
so the graph associated to the correspondence $X$ is
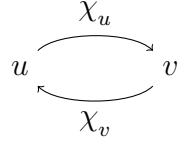
\begin{center}
\begin{tikzpicture}
 \path (-1,0) node (a) {$u$} (1,0) node (b) {$v$};
 \draw[->] (a) .. node[above]{$\Chi_u$} controls (-0.5,0.5) and (0.5,0.5) ..  (b);
 \draw[<-] (a) .. node[below]{$\Chi_v$} controls (-0.5,-0.5) and (0.5,-0.5) ..  (b);
\end{tikzpicture} 
\end{center}
The directed graph $E_\sigma$ associated to $\sigma$ is 
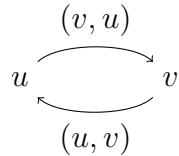
\begin{center}
\begin{tikzpicture}
 \path (-1,0) node (a) {$u$} (1,0) node (b) {$v$};
 \draw[->] (a) .. node[above] {$(v,u)$}  controls (-0.5,0.5) and (0.5,0.5) ..  (b);
 \draw[<-] (a) .. node[below] {$(u,v)$} controls (-0.5,-0.5) and (0.5,-0.5) ..  (b);
\end{tikzpicture} 
\end{center}
\end{ex}

In the proof of Theorem~\ref{characterize} 
we constructed a graph $E$ from a $C^*$-correspondence $X$ 
and showed that $X_E$ recovers $X$, up to isomorphism, from $E$.
Given a graph $F$, 
the same construction recovers $F$, up to isomorphism,
from $X_F$.

\begin{prop}\label{recover-F}
Let $F= (V,F^1,r,s)$ be a directed graph,
and let $X_F$ be the associated graph correspondence.
Then every graph $E$ constructed from $X_F$ 
as in the proof of Theorem~\ref{characterize}
is isomorphic to $F$ 
via a graph isomorphism which fixes the vertex set. 
\end{prop}

\begin{proof}
For ease of notation, set $X=X_F$. 
For $u,v\in V$, $\xi\in X$ and $f\in F^1$, we have
\begin{align*}
(p_u\mdot\xi\mdot p_v)(f) 
&= \Chi_u(r(f))\, \xi(f)\, \Chi_v(s(f))\\
&= \begin{cases}
   \xi(f)& \text{ if $r(f)=u$ and $s(f)=v$}\\
   0& \text{otherwise},
  \end{cases}
\end{align*}
and it follows that 
\[
 X_{uv} = p_u\mdot X\mdot p_v 
\cong \ell^2\bigl(\{ f\in F^1 \mid r(f)=u \text{ and }s(f)=v\}\bigr).
\]
Thus, if $E_{uv}$ is an orthonormal basis for $X_{uv}$,
there exists a bijection
\[
 \phi_{uv}\colon E_{uv}\to \{ f\in F^1 \mid r(f)=u \text{ and }s(f)=v\},
\]
and it is straightforward to check that 
if $E=(V,E^1,r,s)$ is the graph constructed from $X$ 
as in the proof of Theorem~\ref{characterize}
by setting $E^1 = \bigcup_{u,v\in V}E_{uv}$,
then there is a graph isomorphism $\phi\colon E\to F$
such that $\phi = \phi_{uv}$ on $E_{uv}$
and $\phi$ fixes $V$.
\end{proof}

\section{The Graph Correspondence Functor}
\label{functoriality}
 
Let $\GG=\GG(V)$ be the category whose objects are directed graphs with vertex set $V$, and whose morphisms are the injective graph morphisms which are the identity on vertices.
Also let $\CC_0=\CC_0(V)$ denote the category of nondegenerate $C^*$-cor\-re\-spon\-dences over the $C^*$-algebra $c_0(V)$
(we impose the nondegeneracy condition because it is automatically satisfied for graph correspondences). 
The morphisms in $\CC_0$ are just the morphisms of $C^*$-correspondences as in~\eqref{psi}.

For a morphism $\p:E\to F$ in $\GG$, 
consider the linear map
$\psi$ which takes
$\spn\{ \Chi_e \mid e\in E^1 \}\subseteq X_E$
to $X_F$ 
and is defined on generators by
\[
 \psi(\Chi_e) = \Chi_{\phi(e)}.
\]
It is easily verified that $\psi$ satisfies \eqref{psi}
(the assumption that $\phi$ be injective on edges ensures
that $\psi$ is isometric),
and thus extends by continuity to a morphism 
$\psi_\phi$ from $X_E$ to $X_F$ in $\CC_0$.
Then it is routine to verify that the map
$\G:\GG\to \CC_0$ defined on objects by $\G(E) = X_E$
and on morphisms by $\G(\phi) = \psi_\phi$
is a functor.  
We call this the \emph{graph-correspondence functor}. 

Note that \thmref{characterize} says precisely that $\G$ is \emph{essentially surjective}: every object of $\CC_0$ is isomorphic to an object in the image of $\G$.  
Moreover, it is clear from the definition that $\G$ is
\emph{faithful}, that is, $\G$ is injective on each hom-set.  
If it were the case that $\G$ were also \emph{full} 
(surjective on each hom-set) then $\G$ would be a category equivalence between $\GG$ and $\CC_0$.  
However, this is too much to ask:
for one thing, the construction 
of the proof of Theorem~\ref{characterize}  
involves an arbitrary choice of basis
(this is why we have avoided using notation like ``$X\mapsto E_X$''),
although it is trying hard to serve as an adjoint for $\G$. 

Also, it is not hard to see directly that $\G$ can fail to be full,
as in the following example.

\begin{ex}
\label{loopex}
Let $E$ be the graph with a single vertex and a single loop edge, so
that we can identify $c_0(V)$ with the complex numbers. Then $X=\G(E)$ is a
one-dimensional Hilbert space, and a morphism from the correspondence $X$ to itself consists of multiplication by a 
complex number of modulus one.
However, there is only one morphism from the graph $E$ to itself, so there are endomorphisms of $X$ which are not of the form $\G(\p)$ for any endomorphism $\p$ of $E$.
\end{ex}

In spite of the above negative result, 
$\G$ makes a surprisingly close connection between $\GG$ and $\CC_0$.  
For instance, it follows from Theorem~\ref{injective} below that (i)~$\G$ is injective on objects (and hence $\G(\GG)$ is a subcategory of $\CC_0$ which is isomorphic to $\GG$, since $\G$ is faithful), and (ii)~$\G$ \emph{reflects isomorphisms}.  

\begin{thm}
\label{injective}
Let $V$ be a countable set, 
and let $E$ and $F$ be objects of $\GG(V)$.
\begin{enumerate}
\item
If $\G(E)=\G(F)$ in $\CC_0(V)$, then $E=F$ in $\GG(V)$. 
 \item
If $\phi\colon E\to F$ in $\GG(V)$ is such that 
$\G(\phi)\colon \G(E)\to\G(F)$ is an isomorphism in $\CC_0(V)$, 
then $\phi$ is an isomorphism in $\GG(V)$.  
\end{enumerate}
\end{thm}

\begin{proof}
(i)  Suppose $E\neq F$ in $\obj\GG$. 
If there exists $e\in E^1\setminus F^1$, 
then $\Chi_e$ is an element of $\G(E)$ 
which is not in $\G(F)$; 
similarly, if $f\in F^1\setminus E^1$,
then $\Chi_f \in \G(F)\setminus \G(E)$.
Thus if $E^1\neq F^1$, we have $\G(E)\neq \G(F)$.

On the other hand, if $E^1=F^1$, 
then since $E\neq F$
we may choose $e\in E^1=F^1$ 
such that either $r_E(e)\neq r_F(e)$ 
or $s_E(e)\neq s_F(e)$.  
In either case, 
if we set $u=r_E(e)$ and $v=s_E(e)$,
then $p_u\mdot\Chi_e\mdot p_v =\Chi_e$ in $\G(E)$, 
but $p_u\mdot\Chi_e\mdot p_v =0$ in $\G(F)$.  
Thus $\G(E)\neq\G(F)$.

(ii)  For any $\phi\colon E\to F$ in $\GG$
and any $u,v\in V$ we have
\[
 \phi(uE^1v) \subseteq uF^1v,
\]
where we have set $uE^1v = \{ e\in E^1 \mid r(e)=u \text{ and }s(e)=v \}$
and similarly $uF^1v = \{ f\in F^1 \mid r(f)=u \text{ and }s(f)=v \}$.
If we further let
\[
E_{uv} =\{ \Chi_e \mid e\in uE^1v\}\subseteq \G(E)
\quad\text{and}\quad
F_{uv} = \{ \Chi_f \mid f\in uF^1v \}\subseteq \G(F),
\]
it follows that
\begin{align*}
 \G(\phi)(E_{uv})
= \{ \G(\phi)(\Chi_e) \mid e\in uE^1v \}
= \{ \Chi_{\phi(e)} \mid e\in uE^1v \}
\subseteq F_{uv}.
\end{align*}
Now $E_{uv}$ and $F_{uv}$ are orthonormal bases for
the Hilbert spaces $p_u \cdot \G(E)\cdot p_v$ 
and $p_u \cdot \G(F)\cdot p_v$, respectively 
(see the proof of Proposition~\ref{recover-F}).
If $\G(\phi)$ is an isomorphism in $\CC_0$,
then $\G(\phi)$ restricts to a Hilbert space isomorphism 
between these two spaces,
and this forces $\G(\phi)(E_{uv}) = F_{uv}$.
It follows that $\phi(uE^1v) = uF^1v$, and hence
(because morphisms in $\GG$ are injective
and fix vertices by definition)
$\phi$ is an isomorphism in~$\GG$.
\end{proof}

In light of the next theorem,
$\G$ could be called \emph{essentially full}: every morphism in $\CC_0$ 
is ``isomorphic'' in the obvious sense to a morphism in the range of $\G$.
(This usage extends the sense in which the term ``essentially full'' appears in \cite{rosicky}, 
where it was applied to a functor which was surjective on objects.) 

\begin{thm}
\label{close-to-equiv}
Let $V$ be a countable set.
For each morphism $\psi:X\to Y$ in $\CC_0(V)$ 
there exist objects $E$ and $F$
and a morphism $\phi\colon E\to F$ in $\GG(V)$,
and isomorphisms $\upsilon_E$ and $\upsilon_F$ in $\CC_0(V)$,
such that the following diagram commutes in $\CC_0(V)$:
\[
\xymatrix{
X \ar[rr]^-{\upsilon_E}_-\cong\ar[d]_{\psi} 
& &\G(E)\ar[d]^{\G(\phi)}\\
Y \ar[rr]^-{\upsilon_F}_-\cong
& &\G(F).
}
\]
\end{thm}

\begin{proof}
Let $E$ be a graph associated to $X$ as in the proof of Theorem~\ref{characterize}, 
and let $\upsilon_E\colon X\to X_E = \G(E)$ 
be the isomorphism constructed there.
So for each $u,v\in V$, 
we have chosen an orthonormal basis $E_{uv}$ for $X_{uv}$,
and $E^1 = \bigcup_{u,v\in V} E_{uv}$.
Now $\psi(E_{uv})$ is an orthonormal set in $Y_{uv}$ 
since $\psi$ is isometric 
and because $$\psi(p_u\mdot X\mdot p_v) = p_u\mdot \psi(X)\mdot p_v.$$  
So we may extend $\psi(E_{uv})$ to an orthonormal basis $F_{uv}$ of $Y_{uv}$ and set $$F^1 = \bigcup_{u,v\in V}F_{uv}.$$  
As in the proof of \thmref{characterize}, 
this is a disjoint union so we may define $r_F,s_F:F^1\to V$ by  
$$r_F(f)=u, \ s_F(f)=v, \textrm{ if } f \in F_{uv}$$  
to get a graph $F=(V,F^1,r_F,s_F) \in \obj\GG$,
and then the rule $f\mapsto\Chi_f$ 
extends to an isomorphism $\upsilon_F\colon Y\to X_F = \G(F)$ in $\CC_0$.

Now define $\phi:E\to F$ by letting $\phi=\psi$ on $E^1$, 
and letting $\phi$ be the identity on $V$.  
Note that $\phi$ is injective on $E^1$ since $\psi$ is,
so $\phi$ is a morphism in $\GG$. 

To show that the diagram commutes, 
it suffices to show that 
$\G(\phi)\circ\upsilon_E = \upsilon_F\circ\psi$ on $E^1$,
since the linear span of $E^1$ is dense in $X$.  
But for $e\in E^1$ the definitions of the maps involved give
\begin{align*}
 \G(\phi)\circ\upsilon_E(e)
&=\G(\phi)(\Chi_e)
=\Chi_{\phi(e)}
=\Chi_{\psi(e)}
=\upsilon_F(\psi(e))
= \upsilon_F\circ\psi(e).\qedhere
\end{align*}
\end{proof}

\begin{rem}
Combining Theorem~\ref{injective}(ii) with Theorem~\ref{close-to-equiv},
it is not hard to see that $\G$ is injective on isomorphism classes; 
that is, if $\G(E)\cong \G(F)$ in $\CC_0(V)$, 
then $E\cong F$ in $\GG(V)$.
\end{rem}

%



\providecommand{\bysame}{\leavevmode\hbox to3em{\hrulefill}\thinspace}
\providecommand{\MR}{\relax\ifhmode\unskip\space\fi MR }
\providecommand{\MRhref}[2]{%
  \href{http://www.ams.org/mathscinet-getitem?mr=#1}{#2}
}
\providecommand{\href}[2]{#2}

\end{document}